\documentclass[12pt]{amsart}
\usepackage{amsmath,amssymb}
\usepackage{amsfonts}
\usepackage{amsthm}
\usepackage{latexsym}
\usepackage{graphicx}
\usepackage{amssymb,amsmath,txfonts,mathrsfs,mathtools,titletoc, tikz, stmaryrd}
\usepackage[alphabetic]{amsrefs}
\def\p{\partial}

\def\R{\mathbb{R}}

\def\ol{\overline}

\def\vv<#1>{\langle#1\rangle}

\def\XXint#1#2{\setbox0=\hbox{$#1{#2}{\int}$}{#2}\kern-.5\wd0 }

\def\XXint#1#2#3{{\setbox0=\hbox{$#1{#2#3}{\int}$}
     \vcenter{\hbox{$#2#3$}}\kern-.5\wd0}}
\def\pf{{\it Proof:}~}
\newtheorem{theorem}{Theorem}[section]

\newtheorem{definition}[theorem]{Definition}

\def\e{\epsilon}
\def\vv<#1>{\left\langle#1\right\rangle}
\def\ol{\overline}
\newtheorem{thm}{Theorem}[section]

\newtheorem{lem}{Lemma}[section]

\newtheorem{cor}{Corollary}[section]
\theoremstyle{definition}

\theoremstyle{remark}

\newtheorem{rem}{Remark}[section]

\newtheorem{lemma}[theorem]{Lemma}
\numberwithin{equation}{section}
\newtheorem{example}[theorem]{Example}
\begin{document}

\title{The rigidity and stability of  gradient estimates}
\author{Qixuan Hu,  Guoyi Xu$^1$, Chengjie Yu$^2$}
\address{Qixuan Hu\\ Department of Mathematical Sciences\\Tsinghua University, Beijing\\P. R. China}
\email{huqx20@mails.tsinghua.edu.cn}
\address{Guoyi Xu\\ Department of Mathematical Sciences\\Tsinghua University, Beijing\\P. R. China}
\email{guoyixu@tsinghua.edu.cn}
\address{Chengjie Yu\\ Department of Mathematics\\Shantou University\\P. R. China}
\email{cjyu@stu.edu.cn}
\date{\today}

\begin{abstract}
In this note, we obtain the rigidity of the sharp Cheng-Yau gradient estimate for positive harmonic functions on surfaces with nonegative Gaussian curvature, the rigidity of the sharp Li-Yau gradient estimate for  positive solutions to heat equations  and the related estimates for Dirichlet Green's functions on Riemannian manifolds with nonnegative Ricci curvature.  Moreover, we also obtain the corresponding  stability results.
		 		\\[3mm]
\end{abstract}
\thanks{$^1$Research was partially supported by Beijing Natural Science Foundation Z190003, NSFC 11771230 and NSFC 12141103}
\thanks{$^2$Research was partially supported by GDNSF 2021A1515010264 and NSFC 11571215.}
\dedicatory{Dedicated to Professor Peter Li on the occasion of his $70$th birthday}
\maketitle

\titlecontents{section}[0em]{}{\hspace{.5em}}{}{\titlerule*[1pc]{.}\contentspage}
\titlecontents{subsection}[1.5em]{}{\hspace{.5em}}{}{\titlerule*[1pc]{.}\contentspage}
%\contentsline{section}{\numberline{1}something}{}
%\contentsline{subsection}{\numberline{}1.1 other}{}
\tableofcontents

\section{Introduction}
The connection between harmonic functions and geometry can be traced back to the proof of the uniformization theorem of Riemann surfaces.  For manifolds with non-negative Ricci curvature,  one important result relying on harmonic functions, is Cheeger-Gromoll's splitting theorem\cite{CG}.  Its proof reduces to the fact: the existence of a nontrivial harmonic function with constant norm of gradient results in the splitting of  the manifold. This possibly can be viewed as the first geometric application of gradient estimates for harmonic functions.

Yau's Liouville theorem \cite{Yau} inspired a lot of studies in harmonic functions on Riemannian manifolds with non-negative Ricci curvature.  Especially,  Cheng-Yau \cite{CY} proved the local gradient estimate of harmonic functions,  which is not only important as a technical tool, but also provide a universal method to deal with similar analytical objects such as Green's function, heat kernel \cite{LY2}. 

Recall in \cite{Colding},  for the global positive Green's function $G$ on a complete Riemannian manifold with nonnegative Ricci curvature, let $b= G^{\frac{1}{2 -n}}$, we have $|\nabla b|\leq 1$, and $M^n$ is isometric to $\mathbb{R}^n$ if the equality holds at some point in $M^n$. This result also implies: If $G$ and $\ol{G}$ are the global positive Green's functions on $M$ and $\mathbb{R}^n$ with $n\geq 3$, then $\displaystyle G(p,y)\geq \ol{G}(0,\bar{y})$, where $d(p, y)= d(0, \bar{y})$. Moreover if the equality holds at some point, then $M$ is isometric to $\mathbb{R}^n$. Hence, for the global positive Green's function, the sharp estimate and the related rigidity are known.

The philosophy of the above rigidity result is: If the sharp bounds of some estimates for certain kinds of functions (gradient estimates of positive harmonic functions, positive solutions to heat equations or Green's functions) are achieved on some manifold $M^n$, then $M^n$ is isometric to the model space.  In the case of nonnegative curvature, the model space is usually $\mathbb{R}^n$.  

Besides the isometric rigidity, due to Gromov \cite{Gromov}, we can discuss the stability of manifolds (in other words, the Gromov-Hausdorff distance between the manifolds and the model spaces) when the estimate is close to be sharp.

In this note, we consider rigidity and stability  for the Cheng-Yau gradient estimate of positive harmonic functions on surfaces with nonnegative Gaussian curvature,  Li-Yau gradient estimates for positive solutions to heat equations and some estimates of Dirichlet Green's functions on Riemannian manifolds  with $Rc\geq 0$.  For the stability part, we only focus on local stability, more precisely, the stability of geodesic balls in the manifolds. Our results reveal that the geometry of the domain can be rediscovered by the properties of suitable functions defined on it.

The rest of the paper is organized as follows:We first review some preliminary results in Section \ref{sec prelim}, which will be used to establish the rigidity and stability in later sections.
Although the Cheng-Yau gradient estimate for positive harmonic functions is sharp up to decay rate in global sense, it is not sharp pointwisely. The sharp Cheng-Yau gradient on  geodesic balls is only known in $2$-dim case (see \cite{Xu}). In this note, we obtain the related rigidity and stability for geodesic balls on surfaces in Section \ref{sec RS of CY}. It is well-known that Li-Yau's parabolic gradient estimate is sharp as a global inequality on the whole manifold when Ricci curvature is nonnegative. We will also establish the rigidity of Li-Yau's estimate in Section  \ref{sec RS of CY}, which can be viewed as a parabolic analogue of Colding's rigidity result about global Green's functions. Finally, inspired by Colding's rigidity result of global Green's functions, in Section \ref{sec RS of DG}, we study the corresponding rigidity result for Dirichlet Green's functions in geodesic balls,  and obtain the stability result for Dirichlet Green's functions.

\section{Preliminary for local rigidity and stability}\label{sec prelim}
In this section, we recall some preliminary results that will be used in latter sections.

The following lemma is well-known in Riemannian geometry. We present it here for completeness.
\begin{lemma}\label{lem local Lap rigidity}
{Let $(M^n,g)$ be a complete Riemannian manifold with $Rc(g)\geq 0$,  $p\in M$ and $\Omega$ be a star-shaped domain with respect to $p$.  Let  $r(x)= d(p, x)$. If $\displaystyle \Delta r(x)= \frac{n- 1}{r(x)}$ for any $x\in \Omega\setminus \{p\}$ and $\Omega$ contains no cut-point of $p$, then $\Omega$ is isometric to a star-shaped domain $\widehat\Omega$ in $\mathbb{R}^n$. 
}
\end{lemma}
\pf
Applying the Bochner formula to $r$, we get
\begin{equation}\label{eq-Bochner}
\begin{split}
0=\frac12\Delta|\nabla r|^2=|\nabla^2r|^2+\vv<\nabla\Delta r,\nabla r>+Rc(\nabla r,\nabla r)\geq |\nabla^2r|^2-\frac{n-1}{r^2}.
\end{split}
\end{equation}
On the other hand, by the Cauchy-Schwartz inequality, 
\begin{equation}\label{eq-Cauchy}
|\nabla^2r|^2\geq \frac{(\Delta r)^2}{n-1}=\frac{n-1}{r^2}.
\end{equation} 
Combining (\ref{eq-Cauchy}) and (\ref{eq-Bochner}), we know that the equality of (\ref{eq-Cauchy}) holds and hence
\begin{equation}\label{eq-hess-r}
\nabla^2r=\frac{1}{r}(g-dr\otimes dr)
\end{equation}
in $\Omega$. For any normal minimal geodesic $\gamma$ in $\Omega$ starting from $p$, let $E_1=\gamma',E_2,\cdots, E_n$ be a parallel orthonormal frame along $\gamma$. Then, by  (\ref{eq-hess-r})
\begin{equation}
r_{11}=r_{1i}=0\ \mbox{and}\ r_{ij}=\frac1r\delta_{ij}
\end{equation}
for any $i,j\geq 2$. Moreover, 
\begin{equation*}
\begin{split}
0=\frac12\left(|\nabla r|^2\right)_{ij}=r_{ijk}r_k+r_{ki}r_{kj}+R_{iklj}r_kr_l=\frac{dr_{ij}}{dr}+r_{ki}r_{kj}+R_{iklj}r_kr_l=R_{iklj}r_kr_l.
\end{split}
\end{equation*}
Thus, 
$$R(E_i,\gamma',\gamma',E_j)=0$$
for any $i,j=1,2,\cdots,n$. Finally, by Cartan's isometry theorem  (See \cite[Theorem 1.12.8]{Kling}), $\exp_p^{-1}:\Omega\to \widehat\Omega\subset T_pM\backsimeq\R^n$ is an isometry.
\qed

The following theorem follows from \cite{CC-Ann}, which will be used in the rest of the paper.
\begin{thm}\label{thm Hessian equ and flat}
{On a complete Riemannian manifold $(M^n, g)$ with $Rc\geq 0$, let $r(x)= d(p, x)$ and if 
\begin{align}
\fint_{B_1(p)} \left|\nabla^2\left(r^2
\right)- 2g\right|\leq \delta, \nonumber 
\end{align}
then $B_1(p)$ is $\psi(\delta)$-Gromov-Hausdorff close to $B_1(0)$, where $\displaystyle \lim_{\delta\rightarrow 0}\psi(\delta)= 0$. Especially, if $\nabla^2(r^2)- 2g\equiv 0$, then $B_1(p)$ is isometric to $B_1(0)$.
}
\end{thm}

\pf
{It follows from \cite[Proposition $2.80$]{CC-Ann} and \cite[Theorem $3.6$]{CC-Ann}.
}
\qed

\begin{lem}\label{lem Hessian equ and flat}
{If there is a smooth function $f$ defined on a complete Riemannian manifold $(M^n, g)$ with $Rc\geq 0$ such that $\displaystyle \nabla^2 f= c\cdot g$, where $c> 0$ is some constant,  then $(M^n, g)$ is isometric to $\mathbb{R}^n$ and $\displaystyle f(x)= \frac{c}{2}d(p, x)^2+ f(p)$ for some $p\in M^n$.
}
\end{lem}

\pf
{It is the result of \cite[\S 1]{CC-Ann}.
}
\qed

From \cite[Theorem $6.33$]{CC-Ann}, we have the following existence of cut-off functions. 
\begin{lemma}\label{lem suitable cut-off func}
{On a complete Riemannian manifold $(M^n,g)$ with $Rc\geq 0$, there is  a nonnegative cut-off function $\phi$ with $\phi\leq 1$, such that $\phi\big|_{B_{\frac12}(p)}\equiv 1$, $\mathrm{spt}(\phi)\subseteq B_1(p)$, $|\nabla\phi|\leq C$, and $|\Delta\phi|\leq C$, where $C$ is some universal constant.  
}
\end{lemma}

Finally, by reducing to Theorem \ref{thm Hessian equ and flat}, we have the following local stability for Riemannian surfaces in terms of $\Delta r$. % Note that the higher dimension ($n\geq 3$)  case requires $L^2$ estimate of $(\Delta (r^2)- 2n)$, which is the main reason we can not get the stability result in Theorem \ref{thm comp of DHK-stability-2dim} for $n\geq 3$.
\begin{lemma}\label{lem Lap equ and flat-2}
{On a complete Riemannian surface $(M^2, g)$ with $Rc\geq 0$, let $r(x)= d(p, x)$ and if 
\begin{equation}\nonumber 
%\left\{
%\begin{array}{ll}
\fint_{B_1(p)} \left|\Delta\left (r^2\right)- 4\right|\leq \delta, %&n= 2,  \\
%\fint_{B_1(p)} \left|\Delta \left(r^2\right)- 2n\right|^2\leq \delta,&  n \geq 3,
%\end{array} \right.
\end{equation}
then $\displaystyle B_{1}(p)$ is $\psi(\delta)$-Gromov-Hausdorff close to $B_1(0)$, where $\displaystyle \lim_{\delta\rightarrow 0}\psi(\delta)= 0$. 
}
\end{lemma}
\pf
{By Hessian comparison, $2g-\nabla^2\left(r^2\right)$ is nonnegative. So
 $$4- \Delta\left (r^2\right)\geq\left |\nabla^2\left(r^2\right)- 2g\right|.$$
Then, the conclusion follows from Theorem \ref{thm Hessian equ and flat} directly.

}
\qed

\section{Rigidity and stability of Cheng-Yau and Li-Yau estimate}\label{sec RS of CY}
In this section, we prove the rigidity and stability result for Cheng-Yau gradient estimates on surfaces, and the rigidity of Li-Yau gradient estimate.

We first prove the rigidity of the sharp Cheng-Yau gradient for complete Riemannian surfaces with nonnegative Ricci curvature. In fact, we also provide a simpler proof of \cite[Theorem $3.4$]{Xu}.

\begin{thm}\label{thm-new proof and rigidity-2D-Eu}
Let $(M^2,g)$ be a complete Riemannian surface with $Rc\geq 0$ and $u$ be a positive harmonic function on $B_1(p)$. Then the follows hold:
\begin{enumerate}
\item[(1)]. Let $r(x)= d(p,x)$, 
\begin{align}
|\nabla \ln u|(x)\leq \frac{2}{1-r^2(x)}.\label{sharp 2-dim ge}
\end{align}
\item[(2)]. Furthermore,  if there is a point  $x_0\in B_1(p)$ such that 
$$|\nabla \ln u|(x_0)=\frac{2}{1-r^2(x_0)},$$
 then $B_1(p)$ is isometric to the unit ball in $\R^2$ and $u(x)=P(x,y)$ for some $y\in \p B_1(p)$, where $P(x, y)$ is the Poisson kernel of $B_p(1)$.
\item[(3)]. For $\delta> 0$, if
\begin{align}
\fint_{B_1(p)}\left (\frac{4}{(1- r^2)^2}- |\nabla \ln u|^2\right)+  
\ln\left(\frac{2}{(1- r^2)|\nabla \ln u|}\right)\leq \delta, \nonumber 
\end{align}
then $\displaystyle B_{\frac12}(p)$ is $\psi(\delta)$-Gromov-Hausdorff close to $B_{\frac12}(0)$ with $\displaystyle\lim_{\delta\to 0}\psi(\delta)=0$.
\end{enumerate}
\end{thm}

\begin{rem}\label{rem why we get rigidity}
{The inequality (\ref{sharp 2-dim ge}) was proved in \cite{Xu} by a carefully chosen cut-off function.  Here, we firstly give a more direct proof of (\ref{sharp 2-dim ge}) without using cut-off functions, and then use the strong maximum principle to get the rigidity result.
}
\end{rem}
\pf
{\textbf{(1)}.  It was shown in \cite{Xu} that for any positive harmonic function $u$ on $B_1(p)$ in a complete Riemannian surface $(M^2,g)$ with nonnegative Ricci curvature, one has
\begin{equation}
Q\Delta Q\geq 2Q^3+\|\nabla Q\|^2
\end{equation}
where $Q=\|\nabla\ln u\|^2$.  Let $v=\ln Q$. Then,
\begin{equation}\label{eq-v}
\Delta v\geq 2e^v
\end{equation}
on where $Q>0$.  For $\epsilon\in (0, 1 )$, let 
$$v_\e=2\ln\frac{2}{(1- \epsilon)^2-r^2}.$$ 
By Laplacian comparison,
\begin{equation}\label{eq-v-e}
\Delta v_\e\leq 2 e^{v_\e}. 
\end{equation} 
We claim that $v\leq v_\e$ on $B_{1-\e}(p)$. Then, by letting $\e\to 0^+$, we will get (1). If the claim is not true, then there is $x_0\in B_{1-\e}(p)$ such that 
$$v(x_0)-v_\e(x_0)=\max_{x\in B_{1-\e}(p)}(v-v_\e)>0$$ 
by noting that $v$ is upper semi-continuous on $B_1(p) $ and $v_\e|_{\p B_{1-\e}(p)}=+\infty$.  It is also clear that $Q(x_0)>0$ (Otherwise $v(x_0)=-\infty$ which is impossible since $x_0$ is the maximum point of $v-v_\e$). If $x_0$ is not a cut-point of $p$, then, by (\ref{eq-v}) and (\ref{eq-v-e}),
\begin{equation}
0\geq \Delta(v-v_\e)(x_0)\geq  2e^{v(x_0)}-2e^{v_\e(x_0)}>0
\end{equation} 
which is a contradiction.  When $x_0$ is a cut-point of $p$, the classical Calabi trick (see \cite{Ca}) will also give us a contradiction similarly.  This complete the proof of (1).

\textbf{(2)}.  Let $\displaystyle w(x)= 2\ln \frac{2}{1- r(x)^2}$.  If $v=w$ for some point in $B_1(p)$ , let
\begin{equation}
 \Omega=\{x\in B_1(p)\ |\ v(x)=w(x)\}
\end{equation}
It is then clear that $\Omega$ is nonempty and closed in $B_1(p)$ by continuity.  Moreover, for any $x_0\in \Omega$, because $v(x_0)=w(x_0)$, there is an open neighborhood $B_{\delta}(x_0)$ of $x_0$, such that $Q>0$ in $B_{\delta}(x_0)$. Note that
\begin{equation}
0\leq(\Delta v-2e^{v})-(\Delta w-2e^{w})=\Delta(v-w)-c(x)(v-w)
\end{equation}
with 
$$c(x)=\left\{\begin{array}{ll}\frac{2\left(e^{v(x)}-e^{w(x)}\right)}{v(x)-w(x)}& v(x)\neq w(x)\\
2e^{v(x)}&v(x)=w(x).\end{array}\right.$$
in $B_{\delta}(x_0)$. Since $c(x)>0$, by the strong maximum principle (see \cite{Ca} or \cite[Theorem 8.19]{GT}), we know that $B_{\delta}(x_0)\subset\Omega$. This implies that $\Omega$ is open. Now, $\Omega$ as a nonempty open and closed subset of $B_1(p)$ must be $B_1(p)$ since $B_1(p)$ is connected. Hence, $v=w$ all over $B_1(p)$ and moreover the equality of the Laplacian comparison holds. So, $B_1(p)$ is flat and
$$\displaystyle r^2(x)= 1- 2e^{-\frac{v(x)}{2}}$$ 
is smooth in $B_1(p)$ which implies that there is no cut-point of $p$ in $B_1(p)$ (If there is a cut-point $q$ of $p$ in $B_1(p)$, then $q$ is not a conjugate point of $p$ since $B_1(p)$ is flat. So, there are at least two different minimal geodesics joining $p$ to $q$ which contradicts that $\nabla r$ is well defined at $q$ since $r$ is smooth at $q$ ). Now, by Lemma \ref{lem local Lap rigidity}, $B_1(p)$ is isometric to the ball $B_1(0)$ in $\R^2$. Finally by \cite{Xu}, $u$ is the Poisson kernel.

\textbf{ (3)}. It is straightforward to verify that
\begin{equation*}
-\Delta  w+ 2e^{w}\geq\left(1- r^2\right)\left(-\Delta  w+ 2e^{w}\right)= 2\left(4- \Delta\left (r^2\right)\right).
\end{equation*}
So, by (\ref{eq-v}),
\begin{align}
\Delta (v- w)+ 2(e^{w}- e^v)\geq 2\left(4- \Delta \left(r^2\right)\right).  \nonumber 
\end{align}
Let $\phi$ be  the cut-off function in Lemma \ref{lem suitable cut-off func}, and note that $4-\Delta(r^2)\geq 0$ by Laplacian comparison and $w-v\geq 0$ by (1). Then 
\begin{align}
&\quad 2\int_{B_{\frac12}(p)}\left |4- \Delta\left (r^2\right)\right| \nonumber \\
&\leq \int_{B_1(p)}\phi\cdot\left(\Delta (v- w)+ 2(e^{w}- e^v)\right) \nonumber \\
&\leq 2\int_{B_1(p)}\left (\frac{4}{(1- r^2)^2}- |\nabla \ln u|^2\right)- \int_{B_1(p)} 
(w-v)\Delta\phi  \nonumber \\
&\leq 2\int_{B_1(p)} \left(\frac{4}{(1- r^2)^2}- |\nabla \ln u|^2\right)+ 2\sup_{B_1(p)}|\Delta\phi| \int_{B_1(p)} 
\ln\left(\frac{2}{(1- r^2)|\nabla \ln u|}\right) . \nonumber 
\end{align}
Hence,
\begin{align}
&\quad \fint_{B_{\frac12}(p)}\left |4- \Delta\left (r^2\right)\right|\nonumber \\
&\leq C \fint_{B_1(p)}\left (\frac{4}{\left(1- r^2\right)^2}- |\nabla \ln u|^2\right)+  
\ln\left(\frac{2}{(1- r^2)|\nabla \ln u|}\right). \nonumber 
\end{align}
Now the conclusion follows from the above inequality and Lemma \ref{lem Lap equ and flat-2}.
}
\qed

Next, we come to prove the rigidity for the Li-Yau gradient estimate on complete Riemannian manifolds with nonnegative Ricci curvature. We first recall the following Li-Yau's theorem (see \cite{LY2}).
\begin{thm}\label{thm LY gradient est}
{Let $u$ be a positive solution of the heat equation on the  complete Riemannian manifold $(M^n, g)$ with $Rc\geq 0$,  then 
$$\displaystyle (\ln u)_t- |\nabla \ln u|^2+ \frac{n}{2t}\geq 0.$$
}
\end{thm}

The rigidity result  of the Li-Yau gradient estimate on complete Riemannian manifolds with nonnegative Ricci curvature was obtained in \cite{Ni}, we present another proof here. 
\begin{thm}\label{thm strong max prin of gradient est}
{Let $u$ be a positive solution of the heat equation on the complete Riemannian manifold $(M^n, g)$ with $Rc\geq 0$,   and $$\displaystyle (\ln u)_t- |\nabla \ln u|^2+ \frac{n}{2t}= 0$$ at some point $(x_0, t_0)\in M^n\times \mathbb{R}^+$.  Then $(M^n, g)$ is isometric to $\mathbb{R}^n$, and $u(x,t)$ is a multiple of the heat kernel of $\mathbb{R}^n$.
}
\end{thm}

\pf
{Let $f= -\ln u$.  Then 
$$(2t\Delta f- n)(x_0, t_0)= 0.$$  Let $G= t\Delta f- \frac{n}{2}$. It is straightforward to get
\begin{align}
&\left(\frac{\partial}{\partial t}- \Delta\right)(tG)\nonumber \\
=& -\frac{2}{nt^2}(tG)^2- 2\nabla f\cdot \nabla(tG)- 2t\cdot Rc(\nabla f, \nabla f)- \frac{t^2}{n}\left(\sum_{i< j}(f_{ii}- f_{jj})^2\right), \label{1st cru equ}\\
\leq& - 2\nabla f\cdot \nabla(tG). \nonumber
\end{align}
By Theorem \ref{thm LY gradient est}, we get
$$\displaystyle t_0G(x_0, t_0)= t\left(t\Delta f- \frac{n}{2}\right)(x_0, t_0)= \max_{(x, t)\in M^n\times \mathbb{R}^+} tG(x, t).$$ From the strong maximum principle of parabolic PDE (see \cite[Theorem $12.40$]{Chow}),  $tG\equiv 0$ and hence $$\Delta f=\frac{n}{2t}$$  all over $ M^n\times (0,t_0]$. Plugging this into (\ref{1st cru equ}), we have
\begin{align}
Rc(\nabla f, \nabla f)\equiv 0 \quad \quad {\rm and }\quad \quad \nabla^2 f= \frac{1}{2t}g . \nonumber
\end{align}
Combining the above with Lemma \ref{lem Hessian equ and flat}, we get that $M^n$ is $\mathbb{R}^n$, and 
\begin{align}
f(x,t)=\frac{|x-p(t)|^2}{4t}+f(p(t),t), \label{expression of f}
\end{align}
where $p(t)\in \mathbb{R}^n$.  Substituting  (\ref{expression of f}) into the equation of $f$ : $$-f_t-|\nabla f|^2+\frac{n}{2t}=0,$$  we get 
\begin{align}
\frac{n}{2t}- \frac{\partial f}{\partial t}(p(t), t)= - \frac{(x- p(t))\cdot p'(t)}{2t}. \nonumber 
\end{align}
Hence $\frac{d}{dt}p(t)=0$ and $\frac{d}{dt}f(p(t),t)=\frac{n}{2t}$.  Thus $p(t)=p$ is a constant point and $f(x,t)=\frac{|x-p|^2}{4t}+\frac{n}{2}\ln t+C$ for some constant $C\in \mathbb{R}$.
Therefore $u=e^{-f}$ is a multiple of the heat kernel.
}
\qed

As a direct consequence of the rigidity for Li-Yau gradient estimate, we have the following rigidity for Li-Yau's  sharp Harnack inequality of positive solutions to the  heat equation on complete Riemannian manifolds with nonnegative Ricci curvature (See \cite{LY2}).
\begin{cor}\label{cor rigidity of harnack}
Let $u$ be a positive solution of the heat equation on the complete Riemannian manifold $(M^n, g)$ with $Rc\geq 0$. If there exist  $x_1,x_2\in M^n$ and $0<t_1<t_2$, such that 
$$u(x,t_1)=u(x_2,t_2)\left(\frac{t_2}{t_1}\right)^{n/2}\exp\left(\frac{d(x_1, x_2)^2}{4(t_2-t_1)}\right).$$Then $(M^n, g)$ is isometric to $\mathbb{R}^n$ and $u(x,t)$ is a multiple of the heat kernel of $\mathbb{R}^n$.
\end{cor}
\pf Let $\gamma:[t_1,t_2]\to M$ be a minimal geodesic joining $x_1$ to $x_2$. Then, by Theorem \ref{thm LY gradient est} and noting that $|\gamma'|=\frac{d(x_1,x_2)}{t_2-t_1}$, 
\begin{equation*}
\begin{split}
&-\frac{n}{2}\ln\frac{t_2}{t_1}-\frac{d^2(x_1,x_2)}{4(t_2-t_1)}\\
=&\ln u(x_2,t_2)-\ln u(x_1,t_1)\\
=&\int_{t_1}^{t_2}\left(\frac{d}{dt}\ln u(\gamma(t),t)\right)dt\\
=&\int_{t_1}^{t_2}\left((\ln u)_t(\gamma(t),t)+\vv<\nabla\ln u,\gamma'(t)>\right)dt\\
\geq&\int_{t_1}^{t_2}\left(|\nabla\ln u|^2(\gamma(t),t)-\frac{n}{2t}-\frac{d(x_1,x_2)}{t_2-t_1}|\nabla\ln u|(\gamma(t),t)\right)dt\\
=&\int_{t_1}^{t_2}\left(\left(|\nabla\ln u|^2(\gamma(t),t)-\frac{d(x_1,x_2)}{2(t_2-t_1)}\right)^2-\frac{n}{2t}-\frac{d^2(x_1,x_2)}{4(t_2-t_1)^2}\right)dt\\
\geq&-\frac{n}{2}\ln\frac{t_2}{t_1}-\frac{d^2(x_1,x_2)}{4(t_2-t_1)}.
\end{split}
\end{equation*}
Thus, 
$$\displaystyle (\ln u)_t- |\nabla \ln u|^2+ \frac{n}{2t}= 0$$
on $(\gamma(t),t)$ with $t_1\leq t\leq t_2$. Then, by Theorem \ref{thm strong max prin of gradient est},  we complete the proof of the corollary.

\qed

\section{Rigidity and stability for Dirichlet Green's function}\label{sec RS of DG}

For global Green's functions and heat kernels, there are sharp  estimates (see \cite{LY2}, \cite{LTW}, \cite{Colding}).
In this section, we study the sharp estimate and rigidity property of Dirichlet Green's functions, which can viewed as a local version of Colding's rigidity theorem for global positive Green's functions. 

We firstly recall the definition of Dirichlet Green's functions.
\begin{definition}\label{def Dirichlet Green's function}
{On a compact Riemannian manifold $(M^n, g)$ with boundary $\partial M^n\neq \emptyset$, the \textbf{Dirichlet Green's function} $G(x, y)$ is a symmetric function defined on $(M^n\times M^n)\backslash \mathcal{D}(M^n)$ with the properties:
\begin{enumerate}
\item[(i)] for any $f\in C^{\infty}(M^n)$ satisfying $f\big|_{\partial M}= 0$, 
\begin{align}
\int_M G(x, y)\Delta f(y)dy= -f(x) \label{integral equ assump for Green's function}
\end{align}
\item[(ii)] $G(x, y)> 0$ on $(M\backslash \partial M)\times (M\backslash \partial M)$.
\item[(iii)] $G(x, y)= 0$ for any $y\in \partial M^n, x\in M^n\backslash \partial M^n$.
\end{enumerate}
Here $\mathcal{D}(M^n)=\{(x, x): x\in M^n\}\subset M^n\times M^n$.
}
\end{definition}

For Dirichlet Green's functions, we have the following rigidity result related to the gradient of the Dirichlet Green functions, which only applies to higher dimensional case (the dimension is $n\geq 3$).
\begin{lem}\label{lem local rigidity for gradient of DG}
On a complete Riemannian manifold $(M^n, g)$ with $n\geq 3$ and $Rc\geq 0$, let $G(x, y)$ be the Dirichlet Green's function of $B_1(p)$.  Define
$$\displaystyle b(y)= \left(C(n)G(p, y)+ 1\right)^{\frac{1}{2- n}},$$ where $C(n)= (2- n)n\omega_n$.  Then, $$\displaystyle  \sup_{y\in \partial B_1(p)}|\nabla b(y)|\geq \frac{n\omega_n}{V(\partial B_1(p))}.$$
In particular,  if $\displaystyle  \sup_{y\in \partial B_1(p)}|\nabla b(y)|\leq 1$, then $B_1(p)$ is isometric to $B_1(0)\subseteq \mathbb{R}^n$.
\end{lem}
\pf
Note for $y\in \partial B_1(p)$, we have
\begin{align}
\frac{\partial b(y)}{\partial \vec{n}_y}= n\omega_n\frac{\partial G(p, y)}{\partial \vec{n}_y}, \nonumber
\end{align}
where $\partial \vec{n}_y$ is unit inward normal vector at $y$. Recall that $$\displaystyle \int_{\partial B_1(p)}\frac{\partial G(p, y)}{\partial \vec{n}_y}= 1.$$ So
\begin{align}
 \sup_{y\in \partial B_1(p)}|\nabla b(y)|\cdot V(\partial B_1(p))\geq \int_{\partial B_1(p)}|\nabla b|= n\omega_n \int_{\partial B_1(p)}\left|\frac{\partial G(p, y)}{\partial \vec{n}_y}\right|\geq n\omega_n . \nonumber
\end{align}
The first conclusion then follows.

From the rigidity part of the Bishop-Gromov's volume comparison theorem, we get the second conclusion.
\qed

Note that we do not have the global rigidity result as \cite{Colding} for the case $n=2$ because of the following example.
\begin{example}\label{exam 2-dim global rigidity fails}
{Let $M$ be a rotationally symmetric Riemannian surface with center $p$. The metric of $M$ in polar coordinate is $g=dr^2+f(r)^2d\theta^2$.  Then by direct computation, we get the Green function of $M$ at point $p$ has the form $G(p,y)=-\int_1^{r(y)}\frac{1}{f(t)}dt$. If $M$ is Euclidean near $p$, then the Green function of $M$ equals to the Green function of $\mathbb{R}^2$ near $p$, but $M$ may be not isometric to $\mathbb{R}^2$.
}
\end{example}

Now we discuss the sharp $C^0$ estimate of Dirichlet Green's functions and its rigidity and stability result. One difference is that we have the rigidity for all dimension $n\geq 2$ and stability only for $n=2$.

\begin{theorem}\label{thm comp of DHK-stability-2dim}
{On a complete Riemannian manifold $(M^n, g)$ with $Rc\geq 0$, let $G(x, y)$ be the Dirichlet Green's function of $B_1(p)$, and
\begin{equation}\nonumber 
\ol{G}(p, y)= \left\{
\begin{array}{ll}
\frac{d(p,y)^{2-n}-1}{n(n- 2)\omega_n}, & n\geq 3,  \\
-\frac{\ln d(p, y)}{2\pi}, & n =2. 
\end{array} \right.
\end{equation}
Then 
\begin{enumerate}
\item[(1)]. $\displaystyle G(p, y)\geq \ol{G}(p, y)$;
\item[(2)]. if $\displaystyle G(p, y)= \ol{G}(p, y)$ holds at some point $y\in B_1(p)$, then $B_1(p)$ is isometric to $B_1(0)\subseteq \mathbb{R}^n$;
\item[(3)]. If $n=2$ and $\displaystyle \fint_{B_1(p)} \left|G(p, y)-\ol G(p,y) \right|\leq \delta$, then $\displaystyle B_{\frac12}(p)$ is $\psi(\delta)$-Gromov-Hausdorff close to $B_{\frac12}(0)$ with $\displaystyle\lim_{\delta\to 0}\psi(\delta)=0$.
\end{enumerate}
}
\end{theorem}

\begin{rem}\label{rem ellp proof for Diri Green}
{Because $\displaystyle G(x, y)= \int_1^\infty H(x, y, t)dt$, the inequality and the rigidity result for Dirichlet Green's functions can be derived from \cite[Lemma $13.6$]{Li-book}. Here we give an elliptic proof for it.
}
\end{rem}
\pf
\textbf{(1)}.  By Laplacian comparison, we have $\Delta \ol{G}\geq 0$. It is well-known that 
\begin{align}
\lim_{x\rightarrow p}\frac{G(p, x)}{\ol{G}(p,  x)}= 1.  \nonumber 
\end{align}
Hence for any $\alpha>0$, there exist $\delta>0$ sufficiently small, such that $$(1+\alpha)G(p,y)-\ol{G}(y)>0$$ for $y\in\partial B_\delta(p)$.
Moreover, $$(1+\alpha)G(p,y)-\ol{G}(y)=0$$ for $y\in\partial B_1(p)$.  So by maximum principle we have 
$$(1+\alpha)G(p,y)\geq \ol{G}(y)$$ on $B_1(p)\setminus B_\delta(p)$. Letting $\alpha\to 0, \delta\to 0$, we get $G(p,y)\geq \ol{G}(y)$. 

\textbf{(2)}. If the equality of (1) holds at some point, then 
$$G(p,y)=\ol{G}(y)$$ for $y\in B_1(p)\setminus{\{p\}}$ by the strong maximum principle. Thus $B_1(p)$ is isometric to $B_1(0)\subset \mathbb{R}^n$ by Lemma \ref{lem local Lap rigidity}.

\textbf{(3)}. Let $r(x)= d(p,  x)$. For simplicity, we write $G(p,\cdot)$ and $\ol G(p,\cdot)$ as $G$ and $\ol G$.  By direct computation,
\begin{equation*}
2n\omega_n r^n\Delta \ol G=2n-\Delta\left(r^2\right).
\end{equation*}
So, 
\begin{align}
2n- \Delta\left(r^2\right)= 2n\omega_n r^n\Delta \left(\ol{G}-G\right).  \nonumber 
\end{align}
Now let  $\phi$ be the cut-off function  as in Lemma \ref{lem suitable cut-off func}. Then, by the Laplacian comparison, and noting that $\Delta \ol G\geq 0$,
\begin{align}
&\int_{B_{\frac{1}{2}}(p)}\left |2n- \Delta\left(r^2\right)\right|
\leq 2n\omega_n \int_{B_1(p)}r^n\Delta \left(\ol G- G\right)\cdot \phi \nonumber \\
&\leq 2n\omega_n \int_{B_1(p)}\Delta\left (\ol G- G\right)\cdot \phi 
=2n\omega_n \int_{B_1(p)}\left(\ol G-{G}\right)\Delta\phi \nonumber \\
& \leq 2n\omega_n \sup_{B_1(p)}|\Delta\phi| \int_{B_1(p)} \left|\ol{G}- G\right|. \nonumber 
\end{align}
Then, the conclusion follows from Lemma \ref{lem Lap equ and flat-2} and the above estimate.
\qed

\section*{Acknowledgments}
Q. Hu and G. Xu thank the Department of Mathematics, Shantou University for hosting their visit to Shantou, part of the work was done during the visit.

\end{document}